\documentclass{amsart}
\usepackage{amssymb,bm}
\usepackage{amsmath}
\usepackage{amsfonts}
\usepackage{amsthm}
\usepackage{color}
\usepackage{graphics,graphicx,url}
\usepackage{enumerate}
\usepackage{comment}
\usepackage{subcaption}
\usepackage{thmtools}
\usepackage{thm-restate}
\usepackage{centernot}
\usepackage[normalem]{ulem}

\usepackage[utf8]{inputenc}
\usepackage{boondox-cal}

\usepackage[linesnumberedhidden,ruled,noend]{algorithm2e}
\newenvironment{subroutine}[1][htb]
  {
   \begin{algorithm}[#1]%
  }{\end{algorithm}}

\usepackage{breqn} 

\usepackage[a4paper]{geometry}
\usepackage{tikz}
\usetikzlibrary{calc,decorations.pathmorphing,patterns,arrows,decorations.markings,positioning}
\usetikzlibrary{automata,chains,fit,shapes}

\usepackage{graphicx,pdflscape,rotating,float}
\usepackage{enumitem}

\newcommand{\bi}{\begin{itemize}}
\newcommand{\ei}{\end{itemize}}
\newcommand{\be}{\begin{enumerate}}
\newcommand{\ee}{\end{enumerate}}
\newcommand{\bc}{\begin{center}}
\newcommand{\ec}{\end{center}}
\newcommand{\bt}{\begin{tabular}}
\newcommand{\et}{\end{tabular}}
\newcommand{\ba}{\begin{array}}
\newcommand{\ea}{\end{array}}

\newcommand{\calA}{\mathfrak A}
\newcommand{\calB}{\mathfrak B}
\newcommand{\calC}{\mathfrak C}

\newtheorem{theorem}{Theorem}

\newtheorem{lemma}[theorem]{Lemma}
\newtheorem{conjecture}[theorem]{Conjecture}

\newtheorem{proposition}[theorem]{Proposition}

\theoremstyle{definition}
\newtheorem{remark}[theorem]{Remark}
\newtheorem{definition}[theorem]{Definition}

\renewcommand{\geq}{\geqslant} \renewcommand{\leq}{\leqslant}

\newcommand{\PSPACE}{\ensuremath{\mathsf{PSPACE}}}

\newcommand{\DSPACE}{\ensuremath{\mathsf{DSPACE}}}
\newcommand{\NP}{\ensuremath{\mathsf{NP}}}

\usepackage{url}
\usepackage{xcolor}

\newcommand{\N}{\mathbb N}
\newcommand{\Z}{\mathbb Z}
\newcommand{\fclrrs}{finite convergent length-reducing rewriting system}
\newcommand{\icfclrrs}{inverse-closed finite convergent length-reducing rewriting system}
\newcommand{\icf}{icfclrrs}

\newcommand{\nT}{{\mathcal n}_T}
\newcommand{\rT}{{\mathcal r}_T}
\newcommand{\lT}{{\mathcal l}_T}
\newcommand{\bndtp}{bndtp}
\renewcommand{\geq}{\geqslant} \renewcommand{\leq}{\leqslant}  
\newcommand{\Oh}{\mathcal{O}}
\usepackage{tikz}
\usetikzlibrary{positioning,calc}
\usepackage{freetikz}

\usepackage{subcaption}

\title[inverse-closed finite convergent length-reducing rewriting systems]{On groups presented by inverse-closed finite \\convergent length-reducing rewriting systems}

\thanks{Research supported by  Australian Research Council grant DP210100271. The authors are grateful for the support and hospitality of the Sydney Mathematical Research Institute (SMRI) 
during a Domestic Visitor Program visit 
in April 2021.}

\author[M. Elder]{Murray Elder}
\address{School of Mathematical and Physical Sciences, University of Technology Sydney, Ultimo NSW 2007, Australia}
\email{murray.elder@uts.edu.au}
\author[A. Piggott]{Adam Piggott}
\address{Mathematical Sciences Institute, Australian National University, Canberra ACT  2601, Australia}
\email{adam.piggott@anu.edu.au}

\date{\today}



\begin{document}

\begin{abstract}
We show that groups presented by inverse-closed finite convergent length-reducing rewriting systems are characterised by a striking geometric property: their Cayley graphs are geodetic and side-lengths of non-degenerate triangles are uniformly bounded. This leads to a new algebraic result: the group is plain (isomorphic to the free product of finitely many finite groups and copies of $\mathbb Z$) if and only if a certain relation on the set of non-trivial finite-order elements of the group is transitive on a bounded set. We use this to prove that deciding if a group presented by an inverse-closed finite convergent length-reducing rewriting system is not plain is in $\mathsf{NP}$. A ``yes'' answer would disprove a longstanding conjecture of Madlener and Otto from 1987. We also prove that the isomorphism problem for plain groups presented by inverse-closed finite convergent length-reducing rewriting systems is in $\mathsf{PSPACE}$.
\end{abstract}

\maketitle

\section{Introduction}
\label{sec:intro}

A group is {\em plain} if it is isomorphic to a free product of finitely many finite groups and finitely many copies of $\Z$. In the 1980s, the following conjecture was framed in an attempt to understand the algebraic structure of groups presented by  finite convergent length-reducing rewriting systems.

\begin{conjecture}[Madlener and Otto \cite{MadlenerOttoLengthReducing}]\label{conj:MO}
A group $G$ admits presentation by a finite convergent length-reducing rewriting system if and only if $G$ is plain.
\end{conjecture}

Diekert \cite{DiekertLengthReducing} showed that the groups presented by finite convergent length-reducing rewriting systems (the {\em fclrrs groups}) 
form a proper subclass of the virtually-free (and hence hyperbolic) groups. Showing that every fclrrs group is plain, or otherwise, has proved difficult. Special cases where the length of the rewriting rules are restricted have been shown, including: if left-hand sides of rules have length at most two, in 1984 by Avenhaus, Madlener and Otto  \cite{AvenhausMadlenerOtto}; if right-hand sides of rules have length at most one, in 2019 by Eisenberg and the second author  \cite{GilmansConjecture}; if right-hand sides of rules have length at most two and the generating set is inverse-closed, in 2020 by the present authors \cite{EP2020}.

One path to resolving Conjecture \ref{conj:MO} requires identifying properties enjoyed by fclrrs groups that distinguish them among the virtually-free groups. In this paper we identify a striking geometric property characterising the inverse-closed fclrrs groups (the {\em \icf\ groups}). A graph is {\em geodetic} if between any pair of vertices there is a unique shortest path.  
A triangle in a graph is said to be {\em  non-degenerate} if its edges are internally disjoint.  We say that a group $G$  has the {\em $k$-bounded non-degenerate triangle property} with respect to a generating set  $\Sigma$
if $k$ is a universal bound on the side-lengths of non-degenerate geodesic triangles in $\Gamma(G, \Sigma)$, the Cayley graph of $G$ with respect to $\Sigma$. 
For a rewriting system $(\Sigma,T)$, let $\lT=\max\{|\ell|_\Sigma\mid (\ell,r)\in T\}$ and $\rT=\max\{|r|_\Sigma\mid (\ell,r)\in T\}$.

\begin{restatable}[Geometric characterisation]{theoremx}{ThmA}
\label{thmA:ndt}
Let $G$ be a group, let $\Sigma$ be an inverse-closed finite generating set for $G$, and let $k\in\N$.  Then $G$ admits presentation by an inverse-closed finite convergent length-reducing rewriting system $(\Sigma,T)$ with $\rT \leq k$ if and only if the Cayley graph $\Gamma(G, \Sigma)$ is geodetic and has the $k$-bounded non-degenerate triangle property.
\end{restatable}

We define a relation $\sim$ on the set of non-trivial finite-order elements in $G$ such that $a \sim b$ if the product $ab$ has finite order.  It follows easily from the normal form theory of free products that in any plain group the relation $\sim$ is transitive. In general, this property does not distinguish the plain groups among the virtually-free groups; for example, if $H$ is any finite group then $\mathbb{Z} \times H$ is a non-plain virtually-free group in which $\sim$ is transitive.  However, using Bass-Serre Theory, information about centralizers in fclrrs groups, and the geometric constraints imposed by Theorem \ref{thmA:ndt}, we prove that the transitivity of $\sim$ characterizes the plain groups among the fclrrs groups (see Lemma \ref{ProperLoopGroups}). For the icfclrrs groups we can sharpen this to checking the transitivity of $\sim$  on a finite set.

\begin{restatable}[Algebraic characterisation]
{theoremx}{ThmB}
\label{thmB:relation}
If $G$ is a group presented by an inverse-closed finite convergent length-reducing rewriting system $(\Sigma, T)$, then the following are equivalent:
\begin{enumerate}
    \item $G$ is plain;
    \item any nontrivial finite-order element in $G$ is contained in a unique maximal finite subgroup of $G$;
    \item the relation $\sim$ is transitive on the set of non-trivial finite-order elements in $G$;
    \item the relation $\sim$ is transitive on the set of non-trivial finite-order elements in $G$ of geodesic length (with respect to $\Sigma$) at most $11\lT$.
\end{enumerate}
\end{restatable}

The equivalence of conditions (1) and (4) in Theorem~\ref{thmB:relation} allows us to reduce the problem of checking whether or not the group presented by $(\Sigma, T)$ is plain to checking whether or not a finite number of elements  have finite order. 
 This can be done efficiently.

\begin{restatable}[Detecting plainness]{theoremx}{ThmD}
\label{thmD:NP} The following decision problem is in \NP: 
 on input an inverse-closed finite convergent length-reducing  rewriting system $(\Sigma, T)$, is  the group  presented by $(\Sigma, T)$  
 not plain?
 \end{restatable}

A further application of Theorem~\ref{thmB:relation} concerns the complexity of the isomorphism problem within the class of virtually-free groups.
Krsti\'{c} \cite{Krstic} proved that the isomorphism problem for  virtually-free groups is decidable. This was later generalised to all hyperbolic groups by Dahmani and  Guirardel \cite{DG}.  There can be no complexity bound for the isomorphism problem when the inputs are given as arbitrary presentations, since deciding if an arbitrary presentation presents the trivial group is undecidable\footnote{That is, suppose one had a bound on deciding if an arbitrary  presentations for two virtually-free (or hyperbolic) groups. Then  given an arbitrary finite presentation for a group $G$, one can input this together with a presentation for the trivial group 
to the hypothetical algorithm and run until the bound.
If  $G$ was trivial, it is virtually-free, so the algorithm would  return ``yes'' within the bound, and if it doesn't, we can conclude the group is  not trivial.}.  
Recent work of 
S{\'e}nizergues and Wei\ss's \cite{SW} shows that the  isomorphism problem in virtually-free groups is decidable in doubly exponential space if the input is a context-free grammar for the word problem, and in $\PSPACE$ if the input is given in the form of  {\em virtually-free presentations}. A virtually-free presentation of a group $G$ specifies a free group $F$ plus a set of representatives $S$ for $F\setminus G$ together with relations describing pairwise multiplications of elements from $F$ and $S$. 
Using the results in this paper we are able to prove the 
same complexity when the input is an inverse-closed finite convergent length-reducing rewriting system presenting a plain group.

 \begin{restatable}[Isomorphism of plain icfclrrs\ groups]{theoremx}{ThmE}
\label{thmE:PSPACE}
\phantom{.}
 The isomorphism problem for plain groups presented by 
 inverse-closed finite convergent length-reducing rewriting systems is decidable in $\PSPACE$.
 \end{restatable}

A contributing factor in the difficulty of Conjecture~\ref{conj:MO} is a paucity of examples of interesting finite convergent length-reducing rewriting systems that present groups. In 1997, Shapiro asked whether or not the plain groups may be characterized as exactly the groups that admit locally-finite geodetic Cayley graphs \cite[p.286]{Shap-Pascal}. Theorems~\ref{thmA:ndt} and~\ref{thmB:relation} are new tools for considering Shapiro's question.

\subsection*{Acknowledgements}
The authors thank  Volker Diekert and Armin Wei\ss\ for helpful feedback on this work.


\section{Preliminaries}

If $\Sigma$ is an {\em alphabet} (a non-empty finite set), we write $\Sigma^\ast$ for the set of finite words over the alphabet $\Sigma$, and $|u|_\Sigma$ for the length of the word $u\in \Sigma^*$;  the empty word, $\lambda$, is the unique word of length $0$.  If $G$ is a group with generating set $\Sigma$, we write $|g|_{G,\Sigma}$ for the length of a shortest word in $(\Sigma\cup\Sigma^{-1})^*$ which evaluates to $g$. We write: $u=v$ if $u,v\in\Sigma^*$ are identical as words;  $u=_Gv$ if $ u,v\in \Sigma^*$ and $u,v$ evaluate to the same element of $G$; and $u=_Gg$ if $ u\in \Sigma^*,g\in G$ and $u$ evaluates  to $g$.  
We write $e_G$ for the identity element in $G$, and $B_{e_G}(r)$ for the subset of $G$ comprising group elements that may be spelled by a word in $(\Sigma\cup\Sigma^{-1})^*$ of length not exceeding $r$.

A length-reducing rewriting system is a pair  $(\Sigma, T)$ where $\Sigma$ is a non-empty alphabet, and $T$ is a subset of $\Sigma^*\times \Sigma^*$, called a set of {\em rewriting rules},
such that for all $(\ell, r) \in T$ we have that $|\ell|_\Sigma > |r|_\Sigma$. We write $\lT= \max\{|\ell|_\Sigma\mid (\ell,r)\in T\}$, and $\rT=\max\{|r|_\Sigma\mid (\ell,r)\in T\}$.

 The set of rewriting rules determines a relation $\to$ on the set $\Sigma^\ast$ as follows: $a \to b$ if $a=u\ell v$,  $b=urv$ and $(\ell, r) \in T$.  The reflexive and transitive closure of $\to$ is denoted $\overset{\ast}{\to}$. A word $u \in \Sigma^\ast$ is {\em irreducible} if no factor is the left-hand side of any  rewriting rule, and hence $u \overset{\ast}{\to} v$ implies that $u = v$.

The reflexive, transitive and symmetric closure of $\to$ is an equivalence denoted $\overset{\ast}{\leftrightarrow}$.  The operation of concatenation of representatives is well defined on the set of $\overset{\ast}{\leftrightarrow}$-classes, and hence makes a monoid $M = M(\Sigma, T)$.  We say that $M$ is the {\em monoid presented by $(\Sigma, T)$}.  When the equivalence class of every letter (and hence also the equivalence class of every word) has an inverse, the monoid $M$ is a group and we say it is {\em the group presented by $(\Sigma, T)$}.  We say that $(\Sigma, T)$ (or just $\Sigma$) is {\em inverse-closed} if for every $a \in \Sigma$, there exists $b \in \Sigma$ such that $ab \overset{\ast}{\to}\lambda$.  Clearly, $M$ is a group when $\Sigma$ is inverse-closed.

A rewriting system $(\Sigma, T)$ is {\em finite} if $\Sigma$ and $T$ are finite sets, and  {\em terminating} (or {\em noetherian)} if there are no infinite sequences of allowable factor replacements.  It is clear that length-reducing rewriting systems are terminating. A rewriting system is called {\em confluent} if whenever $w\overset{\ast}{\to} x$ and $w\overset{\ast}{\to} y$, there exists $z\in \Sigma^*$ such that $x$ and $y$ both reduce to $z$. 
A rewriting system is called {\em convergent} if it is terminating and confluent.  In some literature, finite convergent length-reducing rewriting systems are called finite \emph{Church-Rosser Thue} systems. Since a finite length-reducing rewriting system is necessarily terminating, the well-known Newman's Lemma \cite[p.69]{SRSBookOtto} gives that checking a finite list of words (corresponding to `critical-pairs') is enough to determine whether or not the rewriting system is convergent.  This can be completed in time that is polynomial in the size of the rewriting system.

If $(\Sigma, T)$ is a finite convergent length-reducing rewriting system, then any element of $M(\Sigma, T)$  is represented by a unique irreducible word $w \in \Sigma^\ast$, and the word $w$ is the unique {\em geodesic} (shortest word) among all representatives of the element. 

We say that $(\Sigma, T)$ is {\em normalized} if for any rule $(\ell, r) \in T$ we have that $r$ is irreducible, every proper subword of $\ell$ is irreducible, and $(\ell, r), (\ell, r') \in T$ implies $r = r'$. We say that  $(\Sigma, T)$  has {\em irreducible letters} if each letter in $\Sigma$ is irreducible.  We note that if $(\Sigma, T)$ is an inverse-closed finite convergent length-reducing rewriting system that is not normalized or contains reducible letters, then there exist subsets $\Sigma' \subseteq \Sigma$ and $T' \subseteq T$ such that  $(\Sigma', T')$ is a normalized inverse-closed finite convergent length-reducing rewriting system with $\Sigma'$ containing only irreducible letters that presents the same group as $(\Sigma, T)$.  Moreover it is easy to compute such $\Sigma'$ and $ T'$. Therefore, without loss of generality, we may assume that any inverse-closed finite convergent length-reducing rewriting system we consider is normalized with irreducible letters.  In such a rewriting system, every rewriting rule is either: $(a b, \lambda)$ for some $a, b \in \Sigma$; or $(u, v)$ for some $u, v \in \Sigma^\ast$ with $|u|_\Sigma = 1+ |v|_\Sigma \geq 2$. In particular, either every rule has the form $(a b, \lambda)$, in which case the group presented is a free product of cyclic groups \cite{Cochet}, or $\lT = \rT+1$.

We define the {\em size} of a rewriting system $(\Sigma, T)$ to be \[\nT=|\Sigma|+\sum_{(\ell,r)\in T}|\ell r|. \] Note that $\rT,\lT\in\Oh(\nT)$.


\section{Geometry of groups presented by rewriting systems}

For a group $G$ and a finite generating set $\Sigma$ (we shall always assume that $\Sigma$ does not contain the identity element $e_G$), the \emph{undirected Cayley graph of $G$ with respect to $\Sigma$} is the locally-finite simple undirected graph $\Gamma = \Gamma(G, \Sigma)$ with vertex set $G$ and in which distinct vertices $g, h \in G$ are adjacent if and only if $g^{-1} h \in \Sigma \cup \Sigma^{-1}$. Each path $v_0, v_1, \dots, v_n$ in $\Gamma$ is labeled by a word $a_1 \dots a_n \in (\Sigma \cup \Sigma^{-1})^\ast$ where $a_i =_G v_{i-1}^{-1} v_i$.  A geodesic path in $\Gamma$ from the identity element $e_G$ to $g$ is labelled by a geodesic word  $u\in (\Sigma \cup \Sigma^{-1})^\ast$ with $|u|_\Sigma=|g|_{G,\Sigma}$.
A simple undirected connected graph is {\em geodetic} if any pair of vertices is joined by a unique geodesic path.
If $\Gamma(G,\Sigma)$ is geodetic and $g \in G$, we will denote the unique geodesic word evaluating to $g$ by $\gamma_g$.

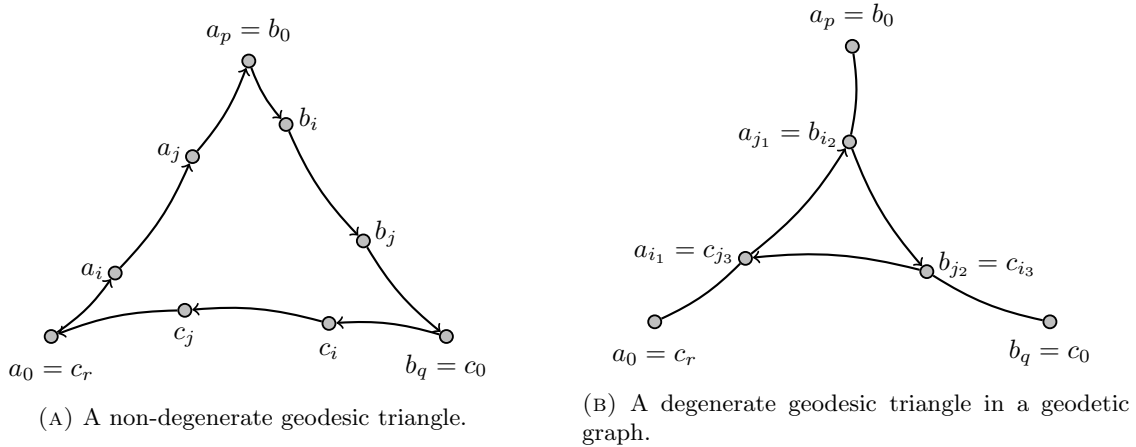
\begin{figure}[!ht]
\begin{subfigure}{.45\textwidth}
  \centering

\begin{tikzpicture}[auto,node distance=.5cm,
    latent/.style={circle,draw,very thick,inner sep=0pt,minimum size=30mm,align=center},
    manifest/.style={rectangle,draw,very thick,inner sep=0pt,minimum width=45mm,minimum height=10mm},
    paths/.style={->, ultra thick, >=stealth'},
]

   \node[dot] [label={[xshift=0cm, yshift=0cm]$a_p=b_0$}](T) {};
             \node[dot] (L) [label={[xshift=0cm, yshift=-.8cm]$a_0=c_r$}, below left= 10em and 7em of T]  {};
             \node[dot] (R) [label={[xshift=0cm, yshift=-.8cm]$b_q=c_0$}, below  right= 10em and 7em of T]  {};
%
     \node[dot] (1i) [label={[xshift=-.3cm, yshift=-.3cm]$a_i$}, above right= 2em and 2em of L]  {};
        \node[dot] (1j) [label={[xshift=-.3cm, yshift=-.3cm]$a_j$}, above right= 4em and 2.5em of 1i]  {};
%
%
      \node[dot] (2i) [label={[xshift=.3cm, yshift=-.3cm]$b_i$}, below right= 2em and 1em of T]  {};
        \node[dot] (2j) [label={[xshift=.3cm, yshift=-.3cm]$b_j$}, below right= 4em and 2.5em of 2i]  {};
%
%
%
      \node[dot] (3i) [label={[xshift=0cm, yshift=-.7cm]$c_i$}, above  left= .1em and 4em of R]  {};
        \node[dot] (3j) [label={[xshift=0cm, yshift=-.7cm]$c_j$},  above left= .1em and 5em  of 3i]  {};

   \path[draw,thick,->]
    (L)  [bend right=10] edge node {} (1i)
                (1i) edge node {} (1j)
                         (1j) edge node {} (T)
    (T) edge node {} (2i)
        (2i) edge node {} (2j)
            (2j) edge node {} (R)
                (R) edge node {} (3i)
                    (3i) edge node {} (3j)
                        (3j) edge node {} (L);

\end{tikzpicture}

\caption{A non-degenerate geodesic triangle.}
  \label{fig:sub-first}
\end{subfigure}
\hspace{5mm}
\begin{subfigure}{.45\textwidth}
  \centering

\begin{tikzpicture}[auto,node distance=.5cm,
    latent/.style={circle,draw,very thick,inner sep=0pt,minimum size=30mm,align=center},
    manifest/.style={rectangle,draw,very thick,inner sep=0pt,minimum width=45mm,minimum height=10mm},
    paths/.style={->, ultra thick, >=stealth'},
]

   \node[dot] [label={[xshift=0cm, yshift=0cm]$a_p=b_0$}](T) {};
             \node[dot] (L) [label={[xshift=0cm, yshift=-.8cm]$a_0=c_r$}, below left= 10em and 7em of T]  {};
             \node[dot] (R) [label={[xshift=0cm, yshift=-.8cm]$b_q=c_0$}, below  right= 10em and 7em of T]  {};

     \node[dot] (1i) [label={[xshift=-.8cm, yshift=-.3cm]$a_{i_1}=c_{j_3}$}, above right= 2em and 3em of L]  {};
        \node[dot] (1j) [label={[xshift=-.8cm, yshift=-.3cm]$a_{j_1}=b_{i_2}$}, above right= 4em and 3.5em of 1i]  {};
                \node[dot] (2j) [label={[xshift=.8cm, yshift=-.3cm]$b_{j_2}=c_{i_3}$}, below right= 4.5em and 2.5em of 1j]  {};

   \path[draw,thick]
    (L) [bend right=10]  edge node {} (1i)
    (1j) edge node {} (T)
    (2j) edge node {} (R);

   \path[draw,thick,->]
                (1i)[bend right=10]  edge node {} (1j)
        (2j) edge node {} (1i)
        (1j) edge node {} (2j);

\end{tikzpicture}

\caption{A degenerate geodesic triangle in a geodetic graph.}
  \label{fig:sub-second}
\end{subfigure}
\caption{Illustrating non-degenerate and degenerate geodesic triangles as in Definition~\ref{defn:nondegen}.}
	\label{fig:nondegen}
\end{figure}

\begin{definition}[Non-degenerate geodesic triangle]\label{defn:nondegen}
Let $\Delta$ be a  
simple undirected 
graph. A {\em geodesic triangle} in $\Delta$ is the union of three geodesic paths $\alpha=a_0,a_1,\dots, a_p$, $\beta=b_0,b_1,\dots, b_q$, $\gamma=c_0, c_1, \dots, ,c_r$ such that $a_p=b_0$, $b_q=c_0$ and $c_r=a_0$. See Figure~\ref{fig:sub-first}. We denote the geodesic triangle by $(\alpha,\beta,\gamma)$.
The geodesic triangle  is {\em non-degenerate} if the vertices $a_i,b_j,c_k$ are all pairwise distinct for $1\leq i\leq p, 1\leq j\leq q, 1\leq k\leq r$. Otherwise we say it is  {\em degenerate}.
\end{definition}
Note that if $\Delta$ is geodetic, then $\Delta$ is degenerate when there exist $i_1,i_2,i_3,j_1,j_2,j_3\in \N$ with $(i_1, i_2, i_3) \neq (0, 0, 0)$ and such that $a_0=c_r,\dots, a_{i_1}=c_{j_3}, a_{j_1}=b_{i_2},\dots, a_p=b_0, b_{j_2}=c_{i_3},\dots, b_q=c_0$
as illustrated in Figure~\ref{fig:sub-second}.

\begin{definition}[Bounded non-degenerate triangle property]\label{def:bndtp}
Let $\Gamma$ be an undirected graph and $k\in\N$. We say $\Gamma$ has the {\em $k$-bounded non-degenerate triangle property ($k$-\bndtp)} if no non-degenerate geodesic triangle in $\Gamma$ has a side-length exceeding $k$.\end{definition}

We make use of  the following  notion from  \cite{EP2020} (the terminology takes its inspiration from \cite{Broomlike}).

 \begin{definition}[$s$-broomlike \cite{EP2020}]\label{defn:broom}
Let $\Delta$ be a geodetic graph and $s$ a positive integer.
We say that $\Delta$ is {\em $s$-broomlike}
if whenever $a_0, \dots, a_n, b$ is a path comprising distinct vertices such that $a_0, \dots, a_n$ is a geodesic but $a_0, \dots, a_n, b$ is not, then the geodesic from $a_0$ to $b$ is \[a_0, \dots, a_{n-p}, b_{n-p+1}, \dots, b_n=b\] for $p \leq s$ and $b_{n-p+1} \neq a_{n-p+1}$.
\end{definition}

  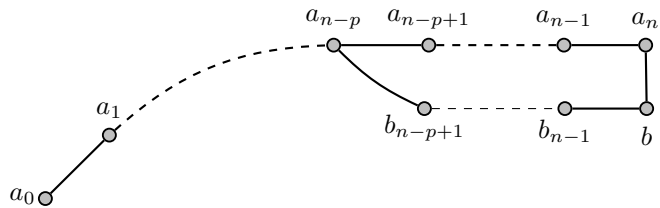
\begin{figure}[!ht]
	\centering

\begin{tikzpicture}[scale=.6,auto,node distance=.5cm,
    latent/.style={circle,draw,very thick,inner sep=0pt,minimum size=30mm,align=center},
    manifest/.style={rectangle,draw,very thick,inner sep=0pt,minimum width=45mm,minimum height=10mm},
    paths/.style={->, ultra thick, >=stealth'},
]

   \node[dot] [label={[xshift=-.3cm, yshift=-.3cm]$a_0$}](0) {};
    \node[dot] (1) [label={[xshift=0cm, yshift=0cm]$a_1$}, above right= 2em and 2em of 0]  {};
        \node[dot] (2) [label={[xshift=0cm, yshift=0cm]$a_{n-p}$},above right = 3em and 8em of 1]  {};
                \node[dot] (3) [label={[xshift=0cm, yshift=-.65cm]$b_{n-p+1}$}, below right = 2em and 3em of 2]  {};
                       
                \node[dot] (5) [label={[xshift=0cm, yshift=0cm]$a_{n-p+1}$},right = 3em of 2]  {};

 \node[dot] (6x) [label={[xshift=0cm, yshift=0cm]$a_{n-1}$}, right = 4.5em of 5]  {};              
 \node[dot] (4x) [label={[xshift=0cm, yshift=-.7cm]$b_{n-1}$},right = 4.7em of 3]  {};

 \node[dot] (6) [label={[xshift=0cm, yshift=0cm]$a_n$}, right = 2.5em of 6x]  {};              
 \node[dot] (4) [label={[xshift=0cm, yshift=-.7cm]$b$},right = 2.5em of 4x]  {};

    \path[draw,thick]
    (0) edge node {} (1)
    
                (6) edge node {} (6x)
                         (2) edge node {} (5)
     (6) edge node {} (4)
              (4) edge node {} (4x);

        \path[draw, thick]
              (2) [bend right=1em] edge node {} (3);

      \draw [dashed] 
     (3) edge node {} (4x);
           
             \draw [thick,dashed] 
              (5) edge node {} (6x)
              
              (1) to [bend left=20]  (2);
      
\end{tikzpicture}

	\caption{Illustrating the $s$-broomlike property (Definition~\ref{defn:broom}).}
	\label{fig:DEFNbroomlike}
\end{figure}

  \begin{lemma}[{\cite[Lemmas~5 and~9]{EP2020}}]
  \label{lem:fclrrs-broom}
If $G$ is presented by an \icfclrrs\ $(\Sigma, T)$, then the undirected Cayley graph of $G$ with respect to $\Sigma$ is geodetic and $\rT$-broomlike.
\end{lemma}

We now prove our first main result.

\ThmA*

\begin{proof}
Let $G$ be presented by an \icf\ $(\Sigma, T)$. 
By Lemma~\ref{lem:fclrrs-broom}, $\Gamma= \Gamma(\Sigma, T)$ is geodetic and $\rT$-broomlike. 
For any path $\epsilon$ in $\Gamma$, we write $|\epsilon|$ for the length of $\epsilon$. Suppose that there exists a non-degenerate geodesic triangle in $\Gamma $ with a side-length exceeding $\rT$.  
Let $f$ be the minimal positive integer such that there exists a non-degenerate triangle $(\alpha_0, \beta_0, \gamma_0)$ in $\Gamma$ such that $|\alpha_0| > \rT$ and $|\alpha_0| + |\beta_0| + |\gamma_0| = f$. Let $\calA$ denote the set of all non-degenerate geodesic triangles $(\alpha, \beta, \gamma)$ in $\Gamma$ such that $|\alpha| > \rT$ and $|\alpha| + |\beta| + |\gamma| = f$.  Let $\calB$ denote the set of all triangles in $\calA$ for which $|\alpha|$ is maximal among all triangles in $\calA$. Since $f$ is a fixed integer, the set $\calB$ is well defined.
Let $\calC$ denote the set of all triangles in $\calB$ for which $|\beta|$ is maximal among all triangles in $B$. Again, since $f$ is fixed, $\calC$ is well defined.

Let $(\alpha, \beta, \gamma) \in \calC$.  Let $\alpha$ be labelled by the word $x_1 x_2 \dots x_p$, and  $\beta$ be labelled by the word $y_1 y_2 \dots y_q$, where $x_i,y_i\in\Sigma$. Without loss of generality we may suppose that the sides are oriented such that 
\[x_1 x_2 \dots x_p y_1 y_2 \dots y_q \gamma
=_G e_G.\] 

First we note that $|\beta|>1$.  If  $|\beta|=1$, then by the  $\rT$-broomlike property $|\alpha|\leq \rT$, which is a contradiction.

The maximality of $|\alpha|$ in $\calA$, gives that $x_1 \dots x_p y_1$ is not reduced---otherwise $(\alpha y_1, y_2 \dots y_q, \gamma)$ would be a triangle in $\calA$ with a longer side.   Since $p > \rT$, the $\rT$-broomlike property gives that there exists $i \geq 1$ and letters $d_{i+1}, \dots d_p \in \Sigma$ such that $x_1 \dots x_i d_{i+1} \dots d_p =_G x_1 x_2 \dots x_p y_1$ and $x_1 \dots x_i d_{i+1} \dots d_p$ is a geodesic.  Let $\epsilon$ be the path from $e_G$ labelled by $x_1 \dots x_i d_{i+1} \dots d_p$.  Since $\gamma^{-1}$ does not start with $x_1$ and $\epsilon$ does, and the paths are geodesics with the same initial point in a geodetic graph, we have that $\gamma$ and $\epsilon$ are internally disjoint. Let $\beta'$ be the geodesic spelled by $y_2 \dots y_q$. See Figure~\ref{fig:new-bndt-proof}.

\begin{figure}[!ht]

  \centering
\begin{tikzpicture}[auto,node distance=.5cm,
    latent/.style={circle,draw,very thick,inner sep=0pt,minimum size=30mm,align=center},
    manifest/.style={rectangle,draw,very thick,inner sep=0pt,minimum width=45mm,minimum height=10mm},
    paths/.style={->, ultra thick, >=stealth'},
]

   \node[dot] (T) {};
             \node[dot] (L) [ below left= 10em and 7em of T]  {};
             \node[dot] (R) [ below  right= 10em and 7em of T]  {};

     \node[dot] (a1) [ above right= 1.5em and 1.5em of L]  {};
     
         \node[dot] (ai) [ below left= 5.8em and 3.1em of T]  {};

         \node[dot] (b1) [ below right= 1.5em and .8em of T]  {};

   \path[draw,thick,->]
    (R)  [bend right=5] edge node {$\gamma$} (L)
    (L)   edge node {$x_1$} (a1)
    (a1) [bend right=10] edge node {} (ai) 
    (ai) edge node {} (T)
    (T) edge node {$y_1$} (b1)
    (b1) edge node {$\beta'$} (R);

   \path[draw,dashed,->]
    (ai) [bend right=13] edge node {} (b1);

\end{tikzpicture}

\caption{Applying the $\rT$-broomlike property to $\alpha y_1$ in the proof of Theorem~\ref{thmA:ndt}.}

	\label{fig:new-bndt-proof}
\end{figure}
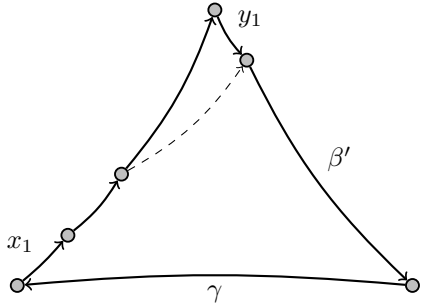

It follows that $(\epsilon, \beta', \gamma)$ is a non-degenerate geodesic triangle with at least one side-length (that of $\epsilon$) exceeding $\rT$ and with $|\epsilon| + |\beta'| + |\gamma| = f-1$.  This contradicts the minimality of $f$. Thus every non-degenerate geodesic triangle in $\Gamma$ has side-lengths at most $\rT$.

Conversely, suppose that $\Gamma$ has the $k$-bounded non-degenerate triangle property.   
Define a set of rewriting rules $T$ as follows:
 \begin{multline*}
T = \{(ab, \lambda) \in \Sigma^\ast \times \Sigma^\ast \mid a, b \in \Sigma \text{ such that } a b=_G e_G\} \cup \\ \{(a_1 a_2 \dots a_n, b_1 b_2 \dots b_{n-1}) \in \Sigma^\ast \times \Sigma^\ast \mid (a_1 a_2 \dots a_{n-1},a_n, b_{n-1}^{-1} b_{n-2}^{-1} \dots b_1^{-1})\\
\text{ labels a non-degenerate geodesic triangle in } \Gamma(G, \Sigma)\}.
 \end{multline*} 
By construction, $(\Sigma, T)$ is a finite, inverse-closed, length-reducing rewriting system with $\rT \leq k$.
Since $\Gamma(G, \Sigma)$ is geodetic, to show that $(\Sigma, T)$ is convergent, it suffices to show that any word labelling a non-geodesic path in $\Gamma(G, \Sigma)$ contains a factor that spells the left-hand side of a rule in $T$.  Let $c = c_1 c_2 \dots c_m \in \Sigma^\ast$ be a word labelling a non-geodesic path in $\Gamma(G, \Sigma)$. Consider first that case that $c$ is not freely reduced.  Then there exists an integer $i$ such that $1 \leq i < m$ and $c_i = c_{i+1}^{-1}$; clearly $(c_i c_{i+1}, \lambda) \in T$. Now consider the case that $c$ is freely reduced.  Let $j$ be the minimal integer such that $c_1 \dots c_j$ is not geodesic; since $\Sigma$ does not contain $e_G$, $j \geq 2$.  Let $i$ be the maximal integer such that $c_{i} \dots c_{j-1}$ is geodesic and $c_i \dots c_j$ is not.  The maximality of $i$ implies that there exists a word $d_i \dots d_{j-1} \in \Sigma^\ast$ with $d_i \neq c_i$ and such that $c_i \dots c_j =_G d_i \dots d_{j-1}$.  It follows that $((c_i,\dots, c_{j-1}), c_j,(d_{j-1}^{-1},\dots, d_i^{-1}))$ labels a non-degenerate geodesic triangle in $\Gamma(G, \Sigma)$.  By hypothesis, $j-i\leq k$.  Thus $(c_i \dots c_j, d_i \dots d_{j-1}) \in T$.

Finally, we establish that $(\Sigma, T)$ presents $G$.  Since $\Sigma$ is inverse-closed, $(\Sigma, T)$ presents a group $\hat{G}$.  By construction,  $\ell =_G r$ for every rule $(\ell, r) \in T$; it follows that $G$ is a quotient of $\hat{G}$.  Any equation that holds in $G$ also holds in $\hat{G}$, because any two words that spell the same element in $G$ will reduce by application of rewriting rules to the unique geodesic representing the group element; it follows that $G \cong \hat{G}$.
 \end{proof}
 
\section{Centralizers of finite order elements}

For a group $G$ and an element $g \in G$,  we write $C_G(g)$ for the centralizer of $g$ in $G$; that is, $C_G(g)=\{t\in G\mid tgt^{-1}=g\}$. In 1988, Madlener and Otto identified the following fact about the centralizers of infinite order elements in fclrrs groups. 

\begin{theorem}[Madlener and Otto  {\cite[Corollary 2.4]{MadlenerOttoCommutativity}}]\label{Mad-centralizer}
Let $G$ be a group presented by a finite convergent length reducing rewriting system  $(\Sigma, T)$.
If $g \in G$ is an element of infinite order, then $C_G(g) \cong \mathbb{Z}$. 
\end{theorem}

Using Theorem~\ref{thmA:ndt}, we can gain further 
information about the centralizers in $G$.  
First we observe the following. For $g,h\in G$, let $\mathfrak C_{g,h}=\{t\in G\mid tgt^{-1}=h\}$.

\begin{lemma}
\label{lem:InfiniteCentralizers}
Let $G$ be a group presented by a finite convergent length-reducing rewriting system and let $g, h$ be non-trivial elements that are conjugate in $G$.  The following are equivalent:
\begin{enumerate}
\item $g$ has infinite order;
\item $C_G(g) \cong \mathbb{Z}$;
\item $C_G(g)$ has infinite order;
\item $\mathfrak C_{g,h}$ has infinite order.
\end{enumerate}
\end{lemma}

\begin{proof} 
That (1) implies (2) is given by Theorem~\ref{Mad-centralizer}. That (2) implies (3) is immediate.  

Next we prove that (3) implies (1).  We prove the contrapositive.  Suppose that $g$ has finite order.  Let $j \in C_G(g)\setminus\{e_G\}$. Then $g \in C_G(j)$.  If $j$ has infinite order, by Theorem \ref{Mad-centralizer} we have that $g$ has infinite order; this contradiction shows that $j$ has finite order and $C_G(g)$ is a torsion subgroup of $G$. In hyperbolic groups, torsion subgroups have finite order \cite[Corollaire 36, Chapitre 8]{FrenchBook}.  Hence $C_G(g)$ has finite order. 

The equivalence of (3) and (4) is an elementary exercise in group theory as follows.
Since $g$ and $h$ are conjugate, $\mathfrak C_{g,h}$ is non-empty; fix $j\in \mathfrak C_{g,h}$. 
For each $t\in G$ we have 
$t \in \mathfrak C_{g,h}$ if and only if $tgt^{-1}=h=jgj^{-1}$ if and only if 
$(j^{-1}t)g(t^{-1}j)=g$  if and only if  $j^{-1}t\in C_G(g)$.  The map $t \mapsto j^{-1} t$ is a bijection from $\mathfrak C_{g,h}$ to  $C_G(g)$, and the equivalence of (3) and (4) follows.
\end{proof}

\begin{lemma}\label{lem:ShortConjugatingElement} Let $G$ be a group presented by an \icfclrrs\ $(\Sigma, T)$.
If $g, h \in B_{e_G}(k)\setminus\{e_G\}$
are conjugate elements of finite order and $t \in G$ is such that $t g t^{-1} =h$, then $|t|_{G,\Sigma} \leq 3\rT+2k$.
\end{lemma}

\begin{proof}
Suppose that $g, h \in B_{e_G}(k)\setminus\{e_G\}$
are conjugate elements of finite order and there exists $t \in G$ such that $t g t^{-1} = h$ and $|t|_{G,\Sigma}\geq 3\rT+2k+1$.  Let $\gamma_t=w, \gamma_g=u, \gamma_h=v$ be the geodesic words for $t,g,h$ respectively, as shown in 
 Figure~\ref{fig:firstLemma2}.

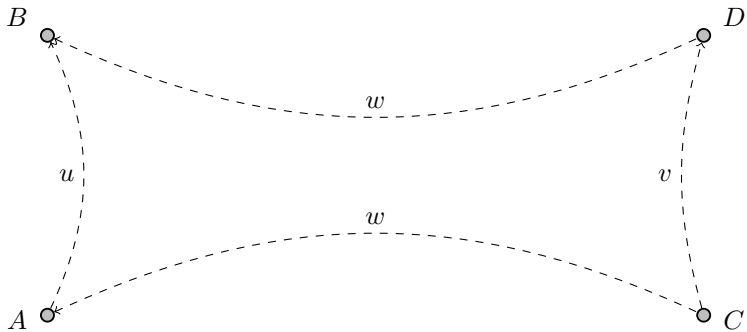
\begin{figure}[!ht]
  \centering
\begin{tikzpicture}[auto,node distance=.5cm,
    latent/.style={circle,draw,very thick,inner sep=0pt,minimum size=30mm,align=center},
    manifest/.style={rectangle,draw,very thick,inner sep=0pt,minimum width=45mm,minimum height=10mm},
    paths/.style={->, ultra thick,  >=stealth'},
]

   \node[dot]  (LB)  [label={[xshift=-.4cm, yshift=-.4cm]$A$}] {};
     \node[dot] (LT) [label={[xshift=-.4cm, yshift=-.1cm]$B$},above= 10em of LB]  {};
     \node[dot] (RB)  [label={[xshift=.4cm, yshift=-.4cm]$C$},right= 24em of LB]  {};
    \node[dot] (RT)  [label={[xshift=.4cm, yshift=-.1cm]$D$},above= 10em of RB]  {};

          \path[draw,dashed]
      (LB) [bend right=25, ->] edge node {$u$} (LT)
        (RB)  [bend left=15] edge node {$v$} (RT)
                (LB) [bend left=25,<-] edge node {$w$} (RB)
    (LT)  [bend right=25] edge node {$w$} (RT);

\end{tikzpicture}
   
\caption{$wuw^{-1}=v$ in $\Gamma(G,\Sigma)$}
	\label{fig:firstLemma2}
\end{figure}

If $x,y$ are labels of vertices in the geodetic graph $\Gamma(G,\Sigma)$, we let $[x,y]$ denote the unique geodesic path starting at $x$ and ending at $y$.

Let $\gamma=[A,D]=\gamma_{t^{-1}h}$ be the geodesic from the vertex marked $A$ to the vertex marked $D$ in Figure~\ref{fig:firstLemma2}. Consider travelling along the two geodesics $\gamma, w^{-1}$ starting at $A$. Let $p$ be the last point visited that lies on both paths. 
Now consider travelling on the two geodesics $w,v$ starting at $C$, and let $p'$ be the last point visited that lies on both paths. Note that $[p',p]$ is a subpath of $w$ since if $p$ was closer to $C$ than $p'$, the geodesic $[p,D]$ includes $p'$ which means it would be shorter to travel from $A$ to $D$ via $[p',D]$, contradicting that $p$ lies on the geodesic from $A$ to $D$.

In this paragraph we prove that $d(p, C) \leq \rT+k$.  We consider cases. Consider first the case that $p=p'$. Then $d(p,C)\leq k$, since $[p',C]$ is a subpath of $v$ which has length at most $k$.  Next consider the case that $p \neq p'$.  The path $[p,p']$ is a side of a non-degenerate geodesic triangle, so has length at most $\rT$ by Theorem~\ref{thmA:ndt}. This means $d(p,C)\leq \rT+k$, since $[p',C]$ is a subpath of $v$ which has length at most $k$.

Now let $\rho=[C,B]=\gamma_{tu}$ be the geodesic from the vertex marked $C$ to the vertex marked $B$ in Figure~\ref{fig:firstLemma2}, and 
let $q$ be the last point visited as you travel along  the two geodesics $\rho, w$ starting at $C$. 
Repeating the above argument we  have  $d(q,A)\leq \rT+k$.

Since $|w|_\Sigma>2(\rT+k)$ we have that $[p,q]$ is a subpath of $w$. Let $\rho_1=[q,B] $, $\gamma_1=[D,p]$, $w_1=[C,p]$, $w_2=[p,q]$,
and $w_3=[q,A]$.
Note that  $\gamma_1w_2\rho_1$ is a geodesic, for if not, it contains a factor of length at most $\rT+1$ which is the left-hand side of a rewrite rule. This factor cannot lie completely inside $\rho$ nor $\gamma$, so does not include the points $p,q$, and this is impossible.
Thus \[w=\gamma_1w_2\rho_1=w_1w_2w_3\] 
See Figure~\ref{fig:secondLemma2}. 

\begin{figure}[!ht]
  \centering
\begin{tikzpicture}[auto,node distance=.5cm,
    latent/.style={circle,draw,very thick,inner sep=0pt,minimum size=30mm,align=center},
    manifest/.style={rectangle,draw,very thick,inner sep=0pt,minimum width=45mm,minimum height=10mm},
    paths/.style={->, ultra thick,  >=stealth'},
]

 \node[dot]  (LB)  [label={[xshift=-.4cm, yshift=-.4cm]$A$}] {};
     \node[dot] (LT) [label={[xshift=-.4cm, yshift=-.1cm]$B$},above= 6em of LB]  {};
     \node[dot] (RB)  [label={[xshift=.4cm, yshift=-.4cm]$C$},right= 26em of LB]  {};
    \node[dot] (RT)  [label={[xshift=.4cm, yshift=-.1cm]$D$},above= 6em of RB]  {};

    \node[dot] (P)  [label={[xshift=0cm, yshift=0cm]$p$}, above left=2em and 6em of RB]  {};
        \node[dot] (Q)  [label={[xshift=0cm, yshift=0cm]$q$}, above right=2em and 6em of LB]  {};

          \path[draw,dashed]
   (LT)  [<-, bend right=10] edge node {$\rho_1$} (Q)
      (LB) [bend right=25] edge node {$u$} (LT)
 (RB)  [bend left=15] edge node {$v$} (RT)
 (LB) [bend left=25,<-] edge node {$w_2$} (RB)
 (RT)  [->, bend left=10, above ] edge node {$\gamma_1$} (P);

       \node (w3)  [above right=0em and 3em of LB]  {$w_3$};
       \node (w1)  [above left=0em and 3em of RB]  {$w_1$};

\end{tikzpicture}
   
\caption{The points $p,q$ and geodesics $\rho=w_1w_2\rho_1, \gamma=w_3^{-1}w_2^{-1}\gamma_1^{-1}$ in $\Gamma(G,\Sigma)$}
	\label{fig:secondLemma2}
\end{figure}
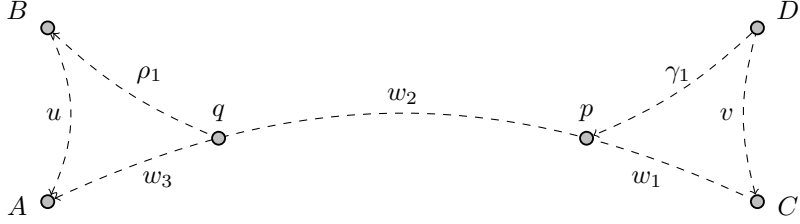

Observe that $|\rho_1|_\Sigma,|\gamma_1|_\Sigma\leq \rT+k$ since they comprise (at most) one side of non-degenerate triangle, plus a factor of $u$ or $v$, and $|u|_\Sigma,|v|_\Sigma\leq k$.

If $|\rho_1|_\Sigma=|w_3|_\Sigma$, then $\rho_1$ and $w_3$ are paths with the same label that share an initial point and so $A = B$; since we assume that $g \neq e_G$, this impossible.  Therefore one of $\rho_1, w_3$ is shorter. Without loss of generality, assume $|w_3|_\Sigma<|\rho_1|_\Sigma$. Then $\rho_1=xw_3$, where $1\leq |x|_\Sigma\leq k+\rT-1$.

Let $y \in \Sigma^\ast$ be the geodesic word such that $w_1 = \gamma_1 y$ (note that $|x|_\Sigma = |y|_\Sigma$).  Then we have $w = \gamma_1 y w_2 w_3 = \gamma_1 w_2 x w_3$; this implies that $y w_2 = w_2 x$.  From this, an exercise in the combinatorics of words 
\cite{har78}
gives that there exist words $r, s \in \Sigma^\ast$ such that $y = rs$, $x = sr$ and $w_2 = (rs)^n r$ for some $n \in \mathbb{N}$.  Thus we have that $w = \gamma_1 (rs)^{n+1}r w_3$.

Let $m$ be a integer such that $m > n$.  We now show that $\gamma_1 (rs)^m r w_3$ is a geodesic word.  Since $|w|_\Sigma \geq 3\rT + 2k +1$ and $|\gamma_1|_\Sigma, |w_3|_\Sigma \leq \rT+k$, we have that $|(rs)^{n+1} r|_\Sigma \geq \rT+1$.  It follows that every length-$(\rT+1)$ factor of $\gamma_1 (rs)^m r w_3$ is also factor of $\gamma_1 (rs)^{n+1} r w_3$.  Since every length-$(\rT+1)$ factor of $\gamma_1 (rs)^{n+1} r w_3$ is a irreducible,  every length every length-$(\rT+1)$ factor of $\gamma_1 (rs)^m r w_3$ is irreducible.  Hence $\gamma_1 (rs)^m r w_3$ is a geodesic word.

Since $w_3 u \rho_1^{-1}$ labels a closed path at $q$, we have that $w_3 u w_3^{-1} r^{-1} s^{-1} =_G e_G$.  It follows that
\[(\gamma_1 (rs)^m r w_3) \ u \ (\gamma_1 (rs)^m r w_3)^{-1} = v\]
for all integers $m > n$.  This implies that $\mathfrak C_{g,h}$ is infinite, which by Lemma~\ref{lem:InfiniteCentralizers} means $g$ does not have finite order, contradicting our hypothesis.
\end{proof}
 
Putting Lemmas~\ref{lem:InfiniteCentralizers} and \ref{lem:ShortConjugatingElement} together we have
\begin{proposition}[Centralizers]
Let $G$ be a group presented by an \icfclrrs\ $(\Sigma, T)$.
Let $g, h \in B_{e_G}(k)\setminus\{e_G\}$ be conjugate elements in $G$.  The following are equivalent:
\begin{enumerate}
\item there exists $t \in G$ such that $t gt^{-1} =h$ and $|t|_{G,\Sigma}>3\rT+2k$;
\item $g$ has infinite order;
\item $C_G(g) \cong \mathbb{Z}$;
\item $C_G(g)$ has infinite order.
\end{enumerate}
\end{proposition}

\begin{remark}
It is known that in a $\delta$-hyperbolic group $G$, if two finite-order elements $g,h \in G$ are conjugate, then there exists an element $x$ conjugating $g$ to $h$ such that the length of $x$ can be bounded in terms of $\delta$, $|\Sigma|$, $|g|_{G,\Sigma}$ and $|h|_{G,\Sigma}$ \cite[Theorem 3.3]{BridsonHowie_Conjugacy} (this result extends to elements of infinite order, and to lists of elements).  The result above bounds the length of \emph{all} elements conjugating $g$ to $h$ in the special case that $G$ is a icfclrrs group.
\end{remark}

\section{Proving Theorem~\ref{thmB:relation}}\label{section:ProvingThmB}

In this section we use the geometric insights of Theorem~\ref{thmA:ndt} and the technical information provided by Lemma~\ref{lem:ShortConjugatingElement} to prove Theorem~\ref{thmB:relation}.  

Recall that $\sim$ is a relation on the set of non-trivial finite-order elements in $G$ defined by the rule  $a \sim b$ if $ab$ has finite order.

\begin{lemma}\label{ProperLoopGroups}
If $G$ is a group presented by a \fclrrs\ $(\Sigma, T)$, then the following are equivalent:
\begin{enumerate}
    \item $G$ is plain;
    \item any nontrivial finite-order element in $G$ is contained in a unique maximal finite subgroup of $G$;
    \item the relation $\sim$ is transitive on the set of non-trivial finite-order elements in $G$.
\end{enumerate}
\end{lemma}

\begin{proof}
Suppose that $G$ is a group presented by a \fclrrs\ $(\Sigma, T)$.  Since $G$ is virtually-free \cite{DiekertLengthReducing}, it is isomorphic to the fundamental group of a finite graph of groups $\mathcal{G}$ with finite vertex groups (and hence also finite edge groups) \cite{KPS}.  Then $\mathcal{G}$ is a connected graph (multiple edges and loops are allowed) in which each vertex is labelled by a group, and each edge is labelled by a group and equipped with two homomorphisms from its label to the vertex group(s) to which it is incident.  Let $V_1, \dots, V_p$ be the finite groups labelling vertices in $\mathcal{G}$, and let $E_1, \dots, E_q$ be the finite groups labelling edges in $\mathcal{G}$. Since $\mathcal{G}$ is connected, $q \geq p-1$.  Without loss of generality, we may assume that $V_i$ is adjacent to $V_{i+1}$ and that the edge incident to $V_i$ and $V_{i+1}$ is labelled $E_i$ for $1 \leq i < p$ (so the edges $E_1, \dots, E_{p-1}$ form a spanning tree in $\mathcal{G}$).   The fundamental group $G_q$ of $\mathcal{G}$ may then be constructed inductively by a sequence of $p-1$ free products with amalgamation followed by a sequence of $q-p+1$ HNN extensions, as follows:
\begin{itemize}
    \item let $G_1 = V_1$;
    \item for $i = 1, 2, \dots, p-1$, let $G_{i+1} = G_i \ast_{E_i} V_{i+1}$ (two copies of $E_i$, one in $V_i$ and one in $V_{i+1}$, are identified);  
    \item for $i = p, p+1, \dots, q$, let $G_{i+1} = G_i \ast _{\phi}$, where $\phi$ is an epimorphism mapping one copy of $E_i$ to another (an infinite-order stable letter $t_k$ conjugates one copy of $E_i$ to another according to the epimorphism $\phi$; the two copies live in the same vertex group if the edge is a loop, or distinct vertex groups otherwise).
\end{itemize}
We note two consequences of the  normal form theory of free products with amalgamation and HNN extensions:
\begin{itemize}
    \item $G_{i+1}$ contains an isomorphic copy of $G_i$ for each $i = 1, 2, \dots, q-1$.
    \item each maximal finite subgroup in $G$ is isomorphic to one of the vertex groups.
\end{itemize}

Since any non-loop edge may be chosen to be in the spanning tree, any non-loop edge with label $E_k$ that is incident to vertices with labels $V_i$ and $V_j$ corresponds to a subgroup of $G$ that is isomorphic to $V_i \ast_{E_k} V_j$, a free product of finite groups amalgamated over a finite subgroup. Every loop with label $E_k$ and incident to a vertex with label $V_i$ corresponds to a subgroup of $G$ that is isomorphic to $V_i \ast_\phi$, an HNN extension of a finite group in which an infinite-order `stable letter' $t_k$ conjugates one embedded copy of $E_k$ to another. It follows from the construction of the fundamental group of $\mathcal{G}$ that $G$ is plain if and only if all edges are labelled by the trivial group.  It is therefore useful to limit the relative sizes of groups labelling edges.

Next we show that, without loss of generality, we may assume that any non-trivial group labelling an edge in $\mathcal{G}$ is strictly smaller than any group labelling an incident vertex. Different arguments are needed for non-loop edges and loops. The arguments given below were first used in \cite{GilmansConjecture}, but are included here for completeness.

In this paragraph we show that, without loss of generality, we may assume that any non-trivial group labelling a non-loop edge in $\mathcal{G}$ is strictly smaller than any group labelling an incident vertex.    Suppose that $\mathcal{G}$ contains a non-loop edge with label $E_k$ that is incident to vertices with labels $V_i$ and $V_j$ and $1 < |E_k| = |V_i| \leq |V_j|$.   The edge corresponds to a subgroup of $G$ that is isomorphic to $V_i \ast_{E_k} V_j$, but $V_i \ast_{E_k} V_j = V_j$. Simply contracting the edge gives a new graph of groups $\mathcal{G}'$ with one less vertex and a fundamental group that is isomorphic to $G$.  We may therefore, without loss of generality, assume that any group labelling a non-loop edge is strictly smaller than the groups labelling the vertices to which the edge is incident. 

Next we show that, because $G$ is presented by a fclrrws group, any non-trivial group labelling a loop in $\mathcal{G}$ must be strictly smaller than the group labelling the incident vertex. Suppose that $\mathcal{G}$ contains a loop labelled $E_k$ that is incident to a vertex with label $V_i$ such that $1 < |E_k| = |V_i|$. Then $G$ contains a subgroup isomorphic to the HNN extension
$V_i \ast_\phi$, where $\phi: V_i \to V_i$ is an automorphism.  Since $V_i$ is finite, $\phi$ has finite order; say $\phi$ has order $f$.  Let $g$ be a non-identity element in $V_i$.  Then $t^f$ is an infinite order element that commutes with $g$, contradicting Theorem \ref{Mad-centralizer}. 

We now know that each edge in the graph of groups $\mathcal{G}$ may be classified as being of one of three types:
\begin{enumerate}
    \item [(E1)] An edge with label $E_k$ such that $|E_k| = 1$
    \item [(E2)] A non-loop edge with label $E_k$ and incident to vertices with labels $V_i$ and $V_j$ such that $1 < |E_k| < |V_i| \leq |V_j|$.
    \item [(E3)] A loop with label $E_k$ that is incident to a vertex with label $V_i$ such that $1 < |E_k| < |V_i|$;
\end{enumerate}

We are now ready to prove the equivalence of conditions (1) and (2).  If $G$ is plain, it follows from the normal form theory of free products that condition (2) holds.  Suppose that $G$ is not plain. Then $\mathcal{G}$ must contain an edge of type (E2) or an edge of type $(E3)$.  Consider first the case that $\mathcal{G}$ contains an edge of type (E2).  Then $G$ contains a subgroup isomorphic to $A \ast_C B$, where $A$ and $B$ are maximal finite subgroups of $G$ and $1 < |C| < |A| \leq |B|$.  Then $A$ and $B$ are distinct maximal finite subgroups in $G$ with a non-trivial intersection, and condition (2) fails.  Now consider the case that $\mathcal{G}$ contains an edge of type (E3).  Then $G$ contains a subgroup isomorphic to $A \ast_\phi$, where $A$ is a maximal finite subgroup of $G$, $C$ is a proper subgroup of $A$, $\phi$ is an epimorphism $\phi: C \to A$, and $t \in G$ is an infinite-order element such $t^{-1} c t = \phi(t)$ for all $c \in C$.   Then $A$ and $t^{-1} A t$ are distinct maximal finite subgroups in $G$ with a non-trivial intersection, and condition (2) fails. 

The equivalence of conditions (2) and (3) is immediate.
\end{proof}

\begin{lemma}\label{lem:HowManyConjClassMFS}
If $G$ is a plain group presented by a \fclrrs\ $(\Sigma, T)$, then the number of conjugacy classes of maximal finite subgroups in $G$ is bounded above by $\nT^2$.
\end{lemma}

\begin{proof}
Suppose that $G$ is a plain group presented by a \fclrrs\ $(\Sigma, T)$.  If $G$ is torsion free, then the trivial subgroup is the only maximal finite subgroup in $G$ and the result is clear. Assume that $G$ is not torsion free.  It follows immediately from Lemma~\ref{ProperLoopGroups}(2) that non-trivial elements contained in representatives of different conjugacy classes of maximal finite subgroups cannot be conjugate. Hence the number of conjugacy class of non-trivial maximal finite subgroups in $G$ is bounded by the number of conjugacy classes of non-trivial finite-order elements. Madlener and Otto \cite{MadOttoOnGroupsHaving} proved that every non-trivial finite-order element in $G$ is conjugate to a proper prefix of some left-hand side of a rewriting rule in $T$.  It follows that the number of conjugacy classes of non-trivial finite-order elements in $G$ is bounded by $\nT^2$. 
\end{proof}

To prove Theorem \ref{thmB:relation}, it suffices to show that if the relation $\sim$ is not transitive on the set of non-trivial finite-order elements in $G$, then a triple of group elements witnessing the failure of transitivity can be found in $B_{e_G}(11\lT)$.  For this we need to develop our understanding of the geodesic structure of finite subgroups in $G$.   

\begin{definition}[Long and short elements]
We say that $g$ is \emph{long} when $|g|_{G,\Sigma} > \rT$, and \emph{short} otherwise.
\end{definition}

\begin{lemma}
\label{lem:LongSubgroups}
Let $G$ be presented by an \icfclrrs\ $(\Sigma,T)$,  and 
let $H$ be a finite subgroup of $G$. If $H$ contains a long element, then there exists a letter $a\in \Sigma$ and geodesic words $h_1, \dots, h_\ell \in \Sigma^\ast$ so that $\{a^{-1} h_1 a, \dots, a^{-1} h_\ell a\}$ is exactly the set of geodesics representing long elements in $H$, and at least half of all elements in $H$ are represented by geodesics of the form $av$ with $v \in \Sigma^\ast$.
\end{lemma}

\begin{proof}
 For each group element $g \in G$, we write $\gamma_g$ for the unique reduced word in $\Sigma^\ast$ such that $|\gamma_g|_\Sigma=|g|_{G,\Sigma}$.  Suppose that there exists a long element $h_0 \in H$.  Then $\gamma_{h_0} = a u b^{-1}$ for $a,b\in H$  and some geodesic word $u\in\Sigma^*$.
Let $A = \{s \in H\mid \gamma_s \text{ starts with } a\}$, 
$B = \{s \in H \mid \gamma_s \text{ starts with } b\}$, $m_A = |A|$ and $m_B = |B|$.  

For each $h \in H \setminus A$, consider the geodesic for $h_0^{-1}h$. 
If 
$\gamma_{h_0^{-1} h}$ does not start with $b$, then
we have $\gamma_s=v_1v_3, \gamma_{h_0^{-1} h}=v_2v_3$ for $|v_1|,|v_2|>0$ since $\gamma_h$  does not start with $a$, and then $(h_0^{-1},v_1,v_2^{-1})$ is a non-degenerate geodesic triangle with $|h_0^{-1}|>2s-2$ (as shown in Figure \ref{fig:finalPicAA}).  This contradiction shows that the first letter of $\gamma_{h_0^{-1} h}$ is $b$.
\begin{figure}[!ht]
	\centering
\begin{tikzpicture}[auto,node distance=.5cm,
    latent/.style={circle,draw,very thick,inner sep=0pt,minimum size=30mm,align=center},
    manifest/.style={rectangle,draw,very thick,inner sep=0pt,minimum width=45mm,minimum height=10mm},
    paths/.style={->, ultra thick, >=stealth'},
]

   \node[dot] (T) {};
             \node[dot] (L) [ below left= 6em and 8em of T]  {};
             \node[dot] (R) [ below  right= 6em and 8em of T]  {};

        \node[dot] (1j) [above right= 3em and 4.5em of L]  {};

      \node[dot] (3i) [ below  left= .1em and 2.5em of R]  {};
        \node[dot] (3j) [below right= .1em and 2.5em  of L]  {};

   \path[draw,thick,->]
    (L)  [bend right=10] edge node {$v_1$} (1j)
     (1j) edge node {$v_3$} (T)

     (3j) edge node [below]{$u$} (3i)
   
     (L) edge node [below] {$a$} (3j)
     [bend left=10] 
     (R)  edge node [below] {$b$}(3i)
      (R) edge node [above] {$v_2$}  (1j);

\end{tikzpicture}
\caption{Proof of Lemma~\ref{lem:LongSubgroups}.}
	\label{fig:finalPicAA}
\end{figure}
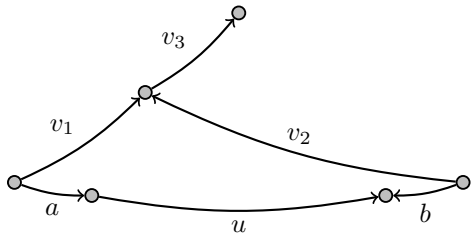

Suppose that $a \neq b$.  Note that distinct elements $j, k  \in H \setminus A$ give distinct geodesics  $\gamma_{h_0^{-1} j}$ and $\gamma_{h_0^{-1} k}$.  It follows that
$m_B \geq |H \setminus A| = |H| - m_A$, so 
\begin{equation}
\label{eqn:twoAA} |H|\leq m_A+m_B.
\end{equation}
Since $e_G \not \in A \cup B$ and $a\neq b$ means $A, B$ are disjoint subsets of $H$, then we have 
\begin{equation*}\label{eqn:one}|H|\leq m_A + m_B < |H|.\end{equation*} 
This contradiction allows us to conclude that $a = b$.

Since $a = b$, we have that for each $h\in H\setminus A$, the first letter of $\gamma_{h_0^{-1} h}$ is $a$. Hence $m_A \geq |H| - m_A$,  so
$m_A\geq  \frac{|H|}2$.  The conclusions of the lemma follow immediately. 
\end{proof}

By Theorem~\ref{thmA:ndt}, if $G$ is presented by a icfclrrs $(\Sigma,T)$, then triangles in $\Gamma(G, \Sigma)$ are $\rT$-thin. A well-known result concerning hyperbolic groups (see \cite{BogGer95, BradyFiniteSubgroups}) then gives that every finite subgroup $H$ of $G$ is conjugate to a subgroup contained within $B_{e_G}(2\rT+1)$.  We can improve this bound for icfclrrs groups.

\begin{proposition}\label{prop:ShortFiniteSubgroups}
Let $G$ be presented by an \icfclrrs\ $(\Sigma,T)$.
Then every finite subgroup $H$ of $G$ is conjugate to a subgroup in $B_{e_G}(\rT+2)$.
\end{proposition}

\begin{proof}Let $H$ be a finite subgroup of $G$.  We shall prove the result by induction on the length of the longest element in $H$.  The result is clearly true in the case that the length of the longest element does not exceed $\rT+2$.  Suppose that the result holds in the case that the length of the longest element in $H$ is at most $n$ for some $n \geq \rT+2$.  Consider the case that the length of the longest element in $H$ is $n+1$.  Let $h_0 \in H$ be such that $|\gamma_{h_0}|_\Sigma=|g|_{G,\Sigma} = n+1$.  Since $n+1 \geq \rT+2 >\rT$, $h_0$ is a long element.  By Lemma~\ref{lem:LongSubgroups}, $\gamma_{h_0} = a u a^{-1}$ for some $a \in \Sigma$ and some $u \in \Sigma^\ast$.  Applying Lemma~\ref{lem:LongSubgroups} to $a^{-1} H a$, because $a^{-1} h_0 a$ is also a long element, yields that $u = b v b^{-1}$ for some $b \in \Sigma$ and some $v \in \Sigma^\ast$.  It follows that $\gamma_{h_0} =  a b v b^{-1} a^{-1}$ and $a \neq b$.  Now we consider the lengths of elements $a^{-1} H a$.  It follows from Lemma~\ref{lem:LongSubgroups} that if $h \in H$ is a long element, then $|\gamma_{a^{-1} h a}| = |\gamma_s|-2$.  It is immediate that if $h \in H$ is a short element (so $|\gamma_h|_\Sigma=|h|_{G,\Sigma} \leq \rT$), then $a^{-1} h a \in B_{e_G}(2s)$. It follows that the length of the longest element in $a^{-1} H a$ does not exceed $\max\{\rT+2, n-1\}$. The inductive hypothesis gives the result.
\end{proof}

\begin{lemma}\label{lem:ElevenLT}
Let $G$ be presented by an \icfclrrs\ $(\Sigma,T)$.
If $J$ is a maximal finite subgroup of $G$ such that $J \cap B_{e_G}(\rT+2) \neq \{e_G\}$, then $J \subseteq B_{e_G}(11 \lT)$.
\end{lemma}

\begin{proof}
If $\rT = 0$, then $G$ is a free group and the statement is vacously true.  We therefore assume that $\rT > 0$.  It follows that $\lT = \rT+1$.  

Suppose that $J$ is a maximal finite subgroup of $G$ such that $J \cap B_{e_G}(\rT+2) \neq \{e_G\}$.  Let $j \in (J \cap B_{e_G}(\rT+2))\setminus \{e_G\}$ and let $k \in J \setminus \{e_G\}$.  We must show that $k \in B_{e_G}(11 \lT)$.   By Proposition~\ref{prop:ShortFiniteSubgroups}, $J$ is conjugate to a subgroup $J'$ contained entirely within $B_{e_G}(\rT+2)$.  Let $t \in G$ be such that $t J' t^{-1} = J$, let $j' = t^{-1} j t$ and $k' = t^{-1} k t$.  By Lemma~\ref{lem:ShortConjugatingElement}, with $g = j'$ and $h = j$, we have that $|t|_{G, \Sigma} \leq 3 \rT + 2(\rT+2)$.  Since $|k'| \leq 2 \rT+2$,  $t^{-1} = t$ and $k = t^{-1} k' t$, we have that \[|k| \leq 2(3 \rT + 2(\rT+2)) + \rT+2 \leq 11 \rT + 10 < 11 (\rT + 1) = 11 \lT.\]
\end{proof}

\ThmB*

\begin{proof}
Suppose that $G$ is a group presented by a \fclrrs\ $(\Sigma, T)$. 
If $\rT=0$, $G$ is the free product of finitely many cyclic groups \cite{Cochet}, but  since $\Sigma$ is inverse-closed, $G$ is free and the theorem is immediate.
So we may assume $\rT>0$ and hence since we have assumed throughout that $(\Sigma, T)$ is normalised, we have $\lT=\rT+1$.

The equivalence of conditions (1), (2) and (3) (in a more general context) was established in Lemma~\ref{ProperLoopGroups}.  We use the notation introduced in the proof of that lemma.  It is clear that condition (3) implies condition (4).  Thus is suffices to show that if (3) does not hold, then (4) does not hold.  

Suppose that $\sim$ is not transitive. We must show that if $\mathcal{G}$ contains an edge of type (E2) or type (E3), then the failure of $\sim$ to be transitive will be witnessed by a triple of elements contained in $B_{e_G}(11\lT)$.

First consider the case that $\mathcal{G}$ contains an edge of type (E2).  It follows that $G$ contains finite subgroups $A, B, C$ with $A$ and $B$ maximal finite subgroups, $A$ not conjugate to $B$, $A \cap B = C$,  $\langle A, B\rangle \cong A \ast_C B$ and $1 < |C| < |A| \leq |B|$. By Proposition~\ref{prop:ShortFiniteSubgroups}, we may assume that $A \subseteq B_{e_G}(\rT+2)$. By Lemma ~\ref{lem:ElevenLT}, we have that $B \subseteq B_{e_G}(11\lT)$.  Let $a \in A \setminus B$, $b \in B \setminus A$ and $c \in C \setminus \{e_G\}$.  Since $\langle A, B\rangle \cong A \ast_C B$, we have that $a \sim c$ and $c \sim b$, but $a \not \sim b$.  Thus the elements $a, b, c$ are in $B_{e_G}(11\lT)$ and witness the failure of $\sim$ to be transitive.

Finally, consider the case that $\mathcal{G}$ contains an edge of type (E3).  It follows that $G$ contains finite subgroups $A, C$ and an infinite order element $t$ with $A$ a maximal finite subgroup, $A \cap t A t^{-1} = C$,  $\langle A, t\rangle \cong A \ast_\phi$ and $1 < |C| < |A|$. By Proposition~\ref{prop:ShortFiniteSubgroups}, we may assume that $A \subseteq B_{e_G}(\rT+2)$.  By Lemma ~\ref{lem:ElevenLT}, we have that $t A t^{-1} \subseteq B_{e_G}(11\lT)$.   Let $a \in A \setminus t A t^{-1}$, $b \in t A t^{-1} \setminus A$ and $c \in A \cap (t A t^{-1}) \setminus \{e_G\}$.  Since $\langle A, t A t^{-1}\rangle \cong A \ast_\phi$, we have that $a \sim c$ and $c \sim b$, but $a \not \sim b$.  Thus the elements $a, b, c$ are in $B_{e_G}(11\lT)$ and witness the failure of $\sim$ to be transitive.
\end{proof}

\section{Algorithms}

We have reduced the problem of deciding plainness to that of deciding whether or not certain elements in the ball of radius $11\lT$ have finite order. 
Narendran and Otto show how to do this in polynomial time.

Recall that we defined the {\em size} of a rewriting system $(\Sigma, T)$ to be \[\nT=|\Sigma|+\sum_{(\ell,r)\in T}|\ell r|_\Sigma. \]

\begin{lemma}
[Theorem 4.8, Narendran and Otto 
\cite{FiniteOrderInPoly}]
\label{lem:finite_order_NarOtto}
Let $G$ be a group presented by an \icfclrrs\ $(\Sigma, T)$ and let $u \in \Sigma^\ast$. There is a deterministic algorithm to decide whether or not $u$ spells an 
element of finite order in $G$, which runs in  time which is polynomial in $|T|, |u|$ and $\mu=\sum_{(r,\ell)\in T} |\ell|_\Sigma$, so polynomial in $|u|\nT$. 
\end{lemma}

\ThmD*

\begin{proof}

If the group $G$ presented by $(\Sigma,T)$ is not plain, we can guess elements 
$u,v,w\in \Sigma^*$ so that $|u|_\Sigma,|v|_\Sigma,|w|_\Sigma\leq 11\lT$, and use Lemma~\ref{lem:finite_order_NarOtto} to verify  that 
\begin{itemize}
    \item each of the words $u, v, w, uv, vw$ spells an element of finite order;
    \item $uw$ spells an element of infinite order.
\end{itemize}
 Since $|u|_\Sigma, |v|_\Sigma, |w|_\Sigma, |uv|_\Sigma, |vw|_\Sigma, |uw|_\Sigma \in\Oh(\nT)$, by Lemma~\ref{lem:finite_order_NarOtto} the algorithm runs in polynomial time in $\nT$.  By Theorem~\ref{thmB:relation}, if the algorithm guesses such a triple, then the relation $\sim$ is not transitive so the group is not plain. If no such triple exists within $(B_{e_G}(11\lT))^3$, then the group is plain so  the algorithm is correct.
\end{proof}

To prove our final theorem, we need to  make use of some simple facts.
\begin{lemma} 
\label{lem:logBound}
Let $(\Sigma, T)$ be an icfclrrs.
Then
$\log_2(|B_{e_G}(\rT+2)|)\leq \nT^2$. 
\end{lemma}
\begin{proof}
Since $\Sigma$ is inverse-closed, there exists a set 
$T'\subseteq T$  comprising exactly one rule $(ax,\lambda)$ with $x\in\Sigma$ for each $a\in\Sigma$.
Then \[\nT\geq |\Sigma|+2|\Sigma|+
\sum_{(\ell,r)\in T\setminus T'} |\ell r|_\Sigma\geq 3+\rT\]
since $|\Sigma|\geq 1$.
We then have that \[|B_{e_G}(\rT+2)|\leq \sum_{i=0}^{\rT+2}|\Sigma|^i\leq |\Sigma|^{\rT+3}\]
so $\log_2(|B_{e_G}(\rT+2)|)\leq (\rT+3)\log_2(|\Sigma|)\leq \nT^2.$
\end{proof}

\begin{lemma}
\label{lem:finiteGenSet}
If $G$ is a finite group, then a minimal generating set for $G$ has  at most $\log_2 |G|$ elements.
\end{lemma}
\begin{proof}
This is a straightforward exercise: let $\{g_1,\dots , g_m\}$ be a minimal generating set, so $g_i\neq e_G$. Let $G_n=\langle g_1,\dots, g_n\rangle$ for $1\leq n\leq m$. Then by minimality $g_{n+1}\not\in G_n$, so there are at least two cosets $e_GG_n, g_{n+1}G_n$, so $|G_{n+1}|\geq 2|G_n|$. By induction $|G|=|G_m|\geq 2^m$.
\end{proof}

We also use this fact.

\begin{lemma}
\label{lem:plainSubgroupMembership}
Let $H_1,\dots, H_k$ be finite groups, 
$ G \cong H_1 \ast H_2 \ast \cdots \ast H_k\ast
           \underbrace{\Z\ast\dots \ast\Z}_{q \text{ copies}},$
and  $g\in H_i$. Then $h\in H_i$ if and only if $h,gh$ have finite order.
\end{lemma}
\begin{proof} This follows immediately from the equivalence of conditions (1) and (2) in  Theorem~\ref{thmB:relation}.
\end{proof}

Recall that   by 
Savitch's theorem \cite{SAVITCH1970177} a problem that can be solved by a nondeterministic algorithm which uses  space $f(n)$  is also in $\DSPACE(f(n)^2)$, so to prove a problem is in \PSPACE\ it suffices to give a nondeterministic polynomial space algorithm.

\begin{proposition}\label{prop:PSPACE-plain}
On input an \icfclrrs\ $(\Sigma, T)$ which presents a plain group $G$, we can output in  space that is polynomial in $\nT$:
\begin{itemize}

\item an integer $k\leq \nT^2$;
\item a list $L_{\Sigma, T}=(v_1,\dots, v_k)$ with each $v_i\in\Sigma^*$ of length at most $\rT+2$;
\item for each  $1\leq i\leq k$, an integer $p_i\leq 2^{\nT^2}$ written in binary; 
\item for each  $1\leq i\leq k$,  a set  $S_i=\left\{u_1,\dots, u_{s_i}\right\}$ with $s_i\leq \nT^2$, $u_j\in \Sigma^*$,  $|u_j|_\Sigma \leq \rT+2$;
 \item an integer $q\leq \nT$;
  \end{itemize}
 so that  for each $1\leq i\leq k$,
$\langle S_i\rangle$ is a maximal finite subgroup $H_i\subseteq B_{e_G}(\rT+2)$ 
of order $p_i$, and 
\begin{gather*}
 G \cong H_1 \ast H_2 \ast \cdots \ast H_k\ast
           \underbrace{\Z\ast\dots \ast\Z}_{q \text{ copies}}.
\end{gather*}
\end{proposition}
\begin{proof}

Assume $\Sigma$ is ordered. 
This induces a shortlex order on  $ \Sigma^*$.

In our subroutines below, when we say ``\textbf{for} $u\in B_{e_G}(a)\setminus B_{e_G}(b)$'' with  $a>b\geq 0$ we mean that the subroutine loops through  in shortlex order every reduced word $u\in \Sigma^*$ with $b<|u|_\Sigma\leq a$. Note that $ B_{e_G}(0)=\{e_G\}$. 
When we say ``check that $u$ has finite/infinite order'' we mean that the subroutine calls the algorithm in 
 Lemma~\ref{lem:finite_order_NarOtto}.

\medskip
First we create an ordered list  $L_{\Sigma,T} = (v_1, v_2, \dots, v_k)$ that will contain exactly one non-trivial element $v_i$ from exactly one representative $H_i\subseteq B_{e_G}(\rT+2)$ of each conjugacy class of maximal finite subgroups in $G$.  By Lemma~\ref{lem:HowManyConjClassMFS}, we know that $k \leq \nT^2$.

\begin{subroutine}
\SetKw{KwSet}{Set} 
    \SetKwInOut{Input}{Input}
    \Input{$(\Sigma, T)$}

 \KwSet{ $L_{\Sigma,T} = ()$.}\\
    \For{$u\in B_{e_G}(\rT + 2) \setminus \{e_G\}$}
      {
      \If{$u$ has finite order} 
       {
    \For{ $t\in B_{e_G}(11 \lT) \setminus B_{e_G}(\rT + 2)$}  {
    Check if $t$ and $tu$ are both finite order.  If both have finite order, by Lemma~\ref{lem:plainSubgroupMembership} $u,t$ lie in the same maximal finite subgroup, so  we are finished considering $u$ (because $u$ lies in a maximal finite subgroup that is not wholly contained in $B_{e_G}(\rT+2)$).  By Lemma~\ref{lem:ElevenLT}, we know that $11\lT$ is enough to check, since the maximal finite subgroup containing $u$ lies completely inside $B_{e_G}(11\lT)$.}
  \For{ $t\in B_{e_G}(5 \rT + 4)\setminus \{e_G\}$} { \For{ $h \in L_{\Sigma,T}$}{
   Check if $tut^{-1}h$ has finite order.  If it does, we are finished considering $u$ (because our list already contains a representative from a maximal finite subgroup in the same conjugacy class as the one containing $u$). Note that by Lemma~\ref{lem:ShortConjugatingElement} (setting $k=\rT+2$) the bound of  $5\rT+4$ is enough to check (since $|u|_\Sigma, |h|_\Sigma\leq \rT+2$).}
 }}
    If $u$ has not been rejected by any of the above steps, append $u$ to $L_{\Sigma,T}$.
    }
    \caption{compute the list  $L_{\Sigma,T}$}
\end{subroutine}

The correctness of this subroutine is guaranteed by  Lemmas~\ref{lem:ElevenLT} and~\ref{lem:ShortConjugatingElement}, and the subroutine runs in \PSPACE.

\medskip
Given a word $v_i$ representing a nontrivial element of a maximal finite subgroup $H_i$ contained wholly in $B_{e_G}(\rT+2)$, we can recover in 
\PSPACE\ the full list of elements of $H_i\setminus\{e_G\}$  as follows:

\begin{subroutine}
    \SetKwInOut{Input}{Input}
    \Input{$(\Sigma, T)$; $v_i\in L_{\Sigma,T}$.}
    \For{$u\in B_{e_G}(\rT + 2) \setminus \{e_G\}$}
      {
     Check that  $u, uv_i$  have finite order.
 If so, return $u$.
    }
    \caption{recover subgroup}
\end{subroutine}
The correctness of this subroutine is guaranteed by   Lemma~\ref{lem:plainSubgroupMembership} and by construction of Subroutine 1.

\medskip
Run Subroutine 1  and store the list $L_{\Sigma, T}$. From this we can read off $k\leq \nT^2$.
Then run Subroutine 2 on each entry $v_i$ in $L_{\Sigma, T}$ with a binary counter to compute the size  $p_i$
of each  $H_i$ (representative of congugacy class of maximal finite subgroup). Since $H_i\subseteq B_{e_G}(\rT+2)$, by Lemma~\ref{lem:logBound} the integers $p_i$ written in binary require space at most $\nT^2$.

\medskip
Now we give another subroutine which on input $v$ in the list  $L_{\Sigma,T}$ can verify in nondeterministic polynomial space  that a given set $S$ is a generating set  for a subgroup $H\subseteq B_{e_G}(\rT+2)$  where $v\in H$.

\newpage

\begin{subroutine}[h!]
    \SetKwInOut{Input}{Input}
    \Input{$v\in L_{\Sigma,T}$; $S=\{u_1,\dots, u_s\}$ where $|u_j|_\Sigma \leq \rT+2$, $s\leq \nT^2$.
    
   (Here suppose $H\subseteq B_{e_G}(\rT+2)$ is the finite subgroup which contains $v$.)} 

\For{$u\in B_{e_G}(\rT + 2) \setminus \{e_G\}$}
{
\If {$u,vu$ have finite order (so $u\in H$)}
{
       \While{$u\neq\lambda$}{
       Nondeterministically choose  $u_i\in S,\epsilon\in\{1,-1\}$.  \\
       Compute the reduced word for $uu_i^\epsilon$ and set $u$ to be this word. 
       
        (Note that if $|u|$ exceeds $\rT+2$ during this procedure, then we have made the wrong guess since we assume $H\subseteq B_{e_G}(\rT+2)$, so we can assume if $S$ is indeed a generating set for $H$, where $H$ is a maximal finite subgroup which lies entirely inside $B_{e_G}(\rT+2)$, that $|u|_\Sigma$ will remain bounded above by $\rT+2$ throughout.)}
}
Once we have verified that $u$ corresponds to an element of $H$ and is equal to a product of letters from $S^{\pm 1}$, we can erase $u$ and move to the next word.
 }
If the subroutine succeeds on every non-empty reduced word $u$ of length at most $\rT+2$, we have verified that $S$ indeed generates $H$.

    \caption{verify generating set}
\end{subroutine}

\medskip

Using this subroutine, we can now compute integers $s_i\leq \nT^2$ and finite generating sets $S_i$ for each $v_i$ in the list $ L_{\Sigma,T}$ using the following nondeterministic polynomial space algorithm.

\begin{subroutine}[h!]
    \SetKwInOut{Input}{Input}
    \Input{$(\Sigma, T)$}
    \For{$1\leq i\leq k$}
      {
       guess an integer $s_i\leq \nT^2$  and a set $S_i=\{u_1, \dots, u_{s_i}\}$ with each $u_j\in \Sigma^*$ reduced and $0<|u_j|_\Sigma\leq \rT+2$.\\
       
      \For{$v_i\in L_{\Sigma,T}$}{
 verify that $S_i$ is a  generating set for the subgroup $H_i$ with $v_i\in H_i$ using Subroutine 3.}
    }
    \caption{compute generating sets}
\end{subroutine}

By Lemma~\ref{lem:finiteGenSet} we are guaranteed that $H_i$ has some generating set $S_i$ of size $s_i\leq \nT^2$, so the subroutine is guaranteed to output a correct answer.

  \medskip
  
Lastly we compute the integer $q$.
Let $G_{\text{ab}}$ denote $G / [G, G]$, the abelianization of $G$.  It is clear that the free-abelian rank of $G_{\text{ab}}$ is equal to $q$, the number of $\mathbb{Z}$ factors in the free product decomposition of $G$.  We may compute the free-abelian rank of the abelianization $G_{\text{ab}}$ from $(\Sigma, T)$ in  time that is polynomial in $\nT$ as follows.  Let $\Sigma' \subseteq \Sigma$ be a subset comprising exactly one generator from inverse pair of inverses.  The information in $(\Sigma, T)$ may be recorded in the form of a group presentation $\langle \Sigma'  \mid R \rangle$, where
$R$ interprets each rewriting rule in $T$ as a relation over the alphabet $(\Sigma')^\pm$. The information in $\langle \Sigma'  \mid R \rangle$ may be encoded in $M$, a $|R| \times |\Sigma'|$ matrix of integers.  These integers record the exponent sums of generators in each relation.  The Smith Normal Form matrix $S$ corresponding to $M$ may be computed in time that is polynomial  in the size of the matrix $S$ \cite{PolynomialSNF}, which is polynomial in $\nT$. The free-abelian rank of $G_{\text{ab}}$ is the number of zero entries along the diagonal of $S$ (see, for example, \cite[pp. 376-377]{SNF}). Note that this means $q\leq \nT$.
\end{proof}
We can now prove:

\ThmE*

\begin{proof}

Two plain groups given as  
\begin{gather*}
  H_1 \ast H_2 \ast \cdots \ast H_k\ast
           \underbrace{\Z\ast\dots \ast\Z}_{q \text{ copies}}, \ \  H_1' \ast H_2' \ast \cdots \ast H_{k'}'\ast
           \underbrace{\Z\ast\dots \ast\Z}_{q' \text{ copies}}
\end{gather*}
are isomorphic if and only if $k=k', q=q'$ and there is a permutation $\sigma\in S^k$ so that $H_i\cong H_{\sigma(i)}'$ for $1\leq i\leq k$. 

Assume $(\Sigma, T), (\Sigma', T')$ are the input, and $N=\max\{\nT,\nT'\}$. 
The procedure we describe  will use polynomial space in $N$, and be nondeterministic.

Guess the following data and store:

\begin{itemize}
   
 \item an integer $q\leq N$; 
\item an integer $k\leq N^2$;
\item a permutation $\sigma$ of length $k$;
\item for each  $1\leq i\leq k$, an integer $p_i\leq 2^{N^2}$ written in binary;
\item for each  $1\leq i\leq k$,  a  set  $S_i=\left\{u_1,\dots, u_{s_i}\right\}$  with $s_i\leq N^2$, $u_j\in \Sigma^*$, and $|u_j|_\Sigma \leq \rT+2$;
\item a list $L'=(z_1,\dots, z_k)$ with $z_i\in (\Sigma')^*$ of length at most $\rT'+2$; 
\item  for each $1\leq i\leq k$,  
 maps $f_i:S_i\to(\Sigma')^*$ with $|f_i(a)|_{\Sigma'}\leq \rT'+2$.
   
\end{itemize}
Note that this requires 
 $\Oh(N^5)$ space: $q\leq N$; $k, |\sigma|\leq N^2$; each $p_i$ in binary requires $N^2$ space and there are $k$ of them, so total $N^4$ space; each set $S_i$ has at most $N^2$ words each of length at most $N$, and there are $k\leq N^2$ such sets so a total of $N^5$ space; $L'$ requires $N^2$ space; each map $f_i$ can be encoded by listing the $s_i\leq N^2$ images of the generators as words of length at most $N$, so $f_i$ requires $N^3$ space and there are $k\leq N^2$ of them so in total $N^5$ space is required.

 Then perform the following tasks.

\begin{enumerate}\item Run the procedure in Proposition~\ref{prop:PSPACE-plain} on input $(\Sigma, T)$  to verify that 
the output includes $k$ generating sets $S_1,\dots, S_k$ with $|\langle S_i\rangle|=p_i$, and the rank of the abelianisation of the plain group presented is $q$.

\item 
Run the procedure in Proposition~\ref{prop:PSPACE-plain} on input $(\Sigma', T')$ to verify that 
the output includes $k$ subgroups having orders $p_{\sigma^{-1}(1)},\dots, p_{\sigma^{-1}(k)}$ (that is, if the $i$-th maximal finite subgroup found by the algorithm is $H_i'$, then $|H_{\sigma(i)}'|=p_i$), the  list $L_{\Sigma',T'}$ is equal to $L'$, and the rank of the abelianisation of the plain group presented is $q$.
\end{enumerate}

So far we have verified that the group presented by $(\Sigma,T)$ is $ H_1 \ast H_2 \ast \cdots \ast H_k\ast
           \underbrace{\Z\ast\dots \ast\Z}_{q \text{ copies}}$ where $|H_i|=p_i$ and $\langle S_i\rangle =H_i$ for $1\leq i\leq k$, and the group presented by $(\Sigma',T')$ is $ H_1' \ast H_2' \ast \cdots \ast H_k'\ast
           \underbrace{\Z\ast\dots \ast\Z}_{q \text{ copies}}$ where $|H_{\sigma(i)}'|=p_i$ and $z_i\in H_i'$ for $1\leq i\leq k$.
           
           To complete the verification that the groups are isomorphic, we need to show $H_i \cong H_{\sigma(i)}'$ for $1\leq i\leq k$.

We do so by showing that each map $f_i$ induces an isomorphism from $H_i$ to $H_{\sigma(i)}'$.
We do this as follows.  Let $u_1\in S_i$ be the first element in the list of generators guessed for $S_i$.
\begin{enumerate}
 

     \item For all pairs $u,v\in \Sigma^*$ of length at most $\rT+2$,  
      \be
     \item[--] if  $u,v$ are reduced, and if  the order of $u, v, uu_1$ and $vu_1$ are finite (using Lemmas~\ref{lem:finite_order_NarOtto} and~\ref{lem:plainSubgroupMembership}),
     \item[--] check that $f_i(uv)=_{G'}f_i(u)f_i(v)$ by writing out both sides as words in $(\Sigma')^*$ and reducing.
     \ee
      This shows that $f_i$ is a homomorphism from $H_i$ to $H_{\sigma(i)}'$.
     \item Check that $f(S_i)$ is a generating set for $H_{\sigma(i)}'$ using Subroutine 3 from 
     the proof of Proposition~\ref{prop:PSPACE-plain} with input $z_{\sigma(i)}$ and $ f(S_i)$.

\end{enumerate}

Once verified, we have that each $f_i$ is a surjective homomorphism between two finite groups of the same size, so $f_i$ is an isomorphism.
\end{proof}

\bibliographystyle{plain}
\bibliography{LengthReducingRewritingBib}

\end{document}